\documentclass{article}

\usepackage{amsmath}
\usepackage{amsthm}

\usepackage[pdftex]{hyperref}
\usepackage{alltt}

\numberwithin{equation}{section}

\renewcommand{\vec}[1]{\mathbf{#1}}

\title{A Method for Fast Diagonalization of\\a 2x2 or 3x3 Real Symmetric Matrix}
\author{M.J. Kronenburg}
\date{}

\begin{document}

\maketitle

\begin{abstract}
A method is presented for fast diagonalization of a 2x2 or 3x3 real symmetric matrix,
that is determination of its eigenvalues and eigenvectors.
The Euler angles of the eigenvectors are computed.
A small computer algebra program is used to compute some of the identities,
and a C++ program for testing the formulas has been uploaded to arXiv.
\end{abstract}


\section{Introduction}

Numerical diagonalization of a real symmetric matrix,
that is determination of its eigenvalues and eigenvectors,
can be found in some textbooks \cite{PTVB07},
and is available in some computer algebra programs \cite{W03}.
This paper shows that in the 2-dimensional and 3-dimensional cases, symbolic diagonalization yields
a fast method to find the eigenvectors.
The formulas for this were published without proof earlier by the author as part of a paper
for an application \cite{K04}.
In this paper emphasis is on the proof of these formulas,
for which a small computer algebra program is used.

\section{Definitions and Basic Identities}

Let $i=1\cdots n$ be the index over the dimensions,
and let $\{\vec{e}_i\}$ be an $n$-dimensional orthonormal basis.
Let $A$ be a square $n$-dimensional real symmetric matrix.
Then there exist a square $n$-dimensional real diagonal matrix $\Lambda$ and
a square $n$-dimensional real orthogonal matrix $D$
such that \cite{PTVB07}:
\begin{equation}\label{diag}
 A = D\cdot\Lambda\cdot D^T
\end{equation}
Let the $\{\lambda_i\}$ be the $n$ real eigenvalues and the $\{\vec{v}_i\}$ be the $n$
mutually orthonormal eigenvectors of $A$,
meaning that:
\begin{equation}
 A\cdot\vec{v}_i = \lambda_i\vec{v}_i
\end{equation}
Then because $A\cdot D=D\cdot\Lambda$, it is clear that $\Lambda$ has as diagonal elements the $\lambda_i$:
\begin{equation}
 \Lambda = \text{diag}(\{\lambda_i\})
\end{equation}
and that the columns of $D$ are the eigenvectors $\vec{v}_i$
with $\vec{v}_i\cdot\vec{v}_j=\delta_{ij}$ \cite{PTVB07}.
The diagonalization problem is that given the symmetric matrix $A$,
how to find the matrices $\Lambda$ and $D$, or in other words
the eigenvalues and the eigenvectors.
The eigenvalues can be found by solving the characteristic equation \cite{PTVB07}:
\begin{equation}
 \text{det}(A-\lambda I) = 0
\end{equation}
which means solving a polynomial equation of degree $n$.
Finding the matrix $D$, that is the eigenvectors, is
more elaborate \cite{PTVB07}.
As the eigenvectors are in pairs $\vec{v}_i$ and $-\vec{v}_i$ and perpendicular,
the orthogonal matrix $D$ can always be a rotation matrix
that rotates the basis vectors $\{\vec{e}_i\}$ onto the eigenvectors $\{\vec{v}_i\}$.
For dimension $n=2$, the 2-dimensional anti-clockwise rotation matrix is:
\begin{equation}\label{rot2}
 D = R(\phi) = \left( 
  \begin{matrix}
    \cos(\phi) & -\sin(\phi) \\
    \sin(\phi) & \cos(\phi) 
  \end{matrix}
 \right)
\end{equation}
For dimension $n=3$, the matrix $D$ is a 3-dimensional rotation matrix,
which is a product of 3 rotation matrices around fixed reference axes:
\begin{equation}\label{rotmat3}
 D = R_1(\phi_1) \cdot R_2(\phi_2) \cdot R_3(\phi_3) 
\end{equation}
where $R_i(\phi_i)$ is the rotation matrix around basis vector $\vec{e}_i$
as the fixed reference axis with the anti-clockwise angle $\phi_i$:
\begin{equation}
 R_1(\phi_1) = \left(
   \begin{matrix}
     1 & 0             & 0            \\
     0 & \cos(\phi_1) & -\sin(\phi_1) \\
     0 & \sin(\phi_1) & \cos(\phi_1)
   \end{matrix}
  \right)
\end{equation}
\begin{equation}
 R_2(\phi_2) = \left(
   \begin{matrix}
     \cos(\phi_2) & 0 & \sin(\phi_2) \\
     0              & 1 & 0             \\
     -\sin(\phi_2) & 0 &  \cos(\phi_2)
   \end{matrix}
  \right)
\end{equation}
\begin{equation}
 R_3(\phi_3) = \left(
   \begin{matrix}
      \cos(\phi_3) & -\sin(\phi_3) & 0 \\
      \sin(\phi_3) & \cos(\phi_3) & 0 \\
     0             & 0            & 1
   \end{matrix}
  \right)
\end{equation}
The method is to symbolically compute $D\cdot\Lambda\cdot D^T$ as functions of the
angles, and then to solve the angles from $A=D\cdot\Lambda\cdot D^T$.
When these angles are known, in the 3-dimensional case,
the Euler rotation angles can be computed,
see section \ref{secteuler}.

\section{The 2-dimensional case}\label{twodimcase}

In the 2-dimensional case the characteristic equation becomes:
\begin{equation}
 \lambda^2 - \lambda(A_{11}+A_{12}) + A_{11}A_{22} - A_{12}^2 = 0
\end{equation}
which yields the following eigenvalues:
\begin{equation}
 \lambda_1 = \frac{1}{2} \left[ A_{11}+A_{22}+\text{sign}(A_{11}-A_{22})
  \sqrt{(A_{11}-A_{22})^2 + 4 A_{12}^2} \right]
\end{equation}
\begin{equation}
 \lambda_2 = \frac{1}{2} \left[ A_{11}+A_{22}-\text{sign}(A_{11}-A_{22})
  \sqrt{(A_{11}-A_{22})^2 + 4 A_{12}^2} \right]
\end{equation}
The matrix $D\cdot\Lambda\cdot D^T$ with $D=R(\phi)$ from (\ref{rot2}) is easily computed:
\begin{equation}\label{mat2}
 D\cdot\Lambda\cdot D^T = \left(
 \begin{matrix}
  \lambda_1\cos(\phi)^2 + \lambda_2\sin(\phi)^2 & (\lambda_1-\lambda_2)\cos(\phi)\sin(\phi) \\
  (\lambda_1-\lambda_2)\cos(\phi)\sin(\phi) & \lambda_1\sin(\phi)^2 + \lambda_2\cos(\phi)^2
 \end{matrix}
 \right)
\end{equation}
and with $2\sin(\phi)\cos(\phi)=\sin(2\phi)$ and
$\arcsin(x)=\arctan(x/\sqrt{1-x^2})$:
\begin{equation}
 \phi = \frac{1}{2}\arctan( \frac{2 A_{12}}{A_{11}-A_{22}} )
\end{equation}
Here $\text{sign}(x)$ is the sign of $x$, and the angle $\phi$ is the anti-clockwise angle
of the eigenvectors with respect to the basis vectors.
The signs are chosen such that interchanging $A_{11}$ and $A_{22}$ leads to interchanging
$\lambda_1$ and $\lambda_2$ and changing the sign of $\phi$.\\
A special case occurs when $A_{11}=A_{22}$ and $A_{12}=0$, in which case $\lambda_1=\lambda_2=\lambda$
and the angle can be taken $\phi=0$.

\section{The 3-dimensional case: the eigenvalues}

In the 3-dimensional case the characteristic equation becomes:
\begin{equation}
 \lambda^3 - b \lambda^2 + c \lambda + d = 0
\end{equation}
where:
\begin{equation}
 b = A_{11}+A_{22}+A_{33}
\end{equation}
\begin{equation}
 c = A_{11}A_{22}+A_{11}A_{33}+A_{22}A_{33}-A_{12}^2-A_{13}^2-A_{23}^2
\end{equation}
\begin{equation}
 d = A_{11}A_{23}^2+A_{22}A_{13}^2+A_{33}A_{12}^2
  - A_{11}A_{22}A_{33} - 2 A_{12}A_{13}A_{23} 
\end{equation}
This cubic polynomial equation is solved by using \cite{Z96}:
\begin{equation}
 p = b^2 - 3c
\end{equation}
and
\begin{equation}
 q = 2b^3 - 9bc - 27d
\end{equation}
When $q^2< 4p^3$ there are three unequal real roots, and when $q^2=4p^3$,
there are three real roots of which at least two are equal.
\begin{equation}
\begin{split}
 p = & \frac{1}{2} \left[ (A_{11}-A_{22})^2 + (A_{11}-A_{33})^2 + (A_{22}-A_{33})^2 \right] \\
 & + 3 ( A_{12}^2 + A_{13}^2 + A_{23}^2 )
\end{split}
\end{equation}
\begin{equation}
\begin{split}
 q = & 18( A_{11}A_{22}A_{33} + 3 A_{12}A_{13}A_{23} )
    + 2( A_{11}^3 + A_{22}^3 + A_{33}^3 ) \\
 & + 9 ( A_{11} + A_{22} + A_{33} )( A_{12}^2 + A_{13}^2 + A_{23}^2 ) \\
 & - 3 ( A_{11} + A_{22} )( A_{11} + A_{33} )( A_{22} + A_{33} ) \\
 & - 27 ( A_{11}A_{23}^2 + A_{22}A_{13}^2 + A_{33}A_{12}^2 )
\end{split}
\end{equation}
\begin{equation}
 \Delta = \arccos( \frac{q}{2\sqrt{p^3}} )
\end{equation}
\begin{equation}
 \lambda_1 = \frac{1}{3} \left[ b + 2\sqrt{p} \cos(\frac{\Delta}{3}) \right]
\end{equation}
\begin{equation}
 \lambda_2 = \frac{1}{3} \left[ b + 2\sqrt{p} \cos(\frac{\Delta+2\pi}{3}) \right]
\end{equation}
\begin{equation}
 \lambda_3 = \frac{1}{3} \left[ b + 2\sqrt{p} \cos(\frac{\Delta-2\pi}{3}) \right]
\end{equation}
Changing the sign of $q$ leads to changing the first plus signs in the $\lambda_i$
into minus signs, which is easily shown by using $\arccos(-x)=\pi-\arccos(x)$ and
$\cos(x)=-\cos(\pi-x)$.

\section{The 3-dimensional case: the eigenvectors}

The eigenvectors are the columns of $D$, which is the product of rotation matrices (\ref{rotmat3}),
which means that when the three angles $\phi_i$ are known, the eigenvectors are known.
These angles are computed by symbollically evaluating the equation
$A=D\cdot\Lambda\cdot D^T$ and solving the angles.
The following two equations are computed in section \ref{program} below by the computer algebra program.
\begin{equation}
 v = \cos(\phi_2)^2 = \frac{A_{12}^2+A_{13}^2+(A_{11}-\lambda_3)(A_{11}+\lambda_3-\lambda_1-\lambda_2)}
    {(\lambda_2-\lambda_3)(\lambda_3-\lambda_1)}
\end{equation}
\begin{equation}
 w = \cos(\phi_3)^2 = \frac{A_{11}-\lambda_3+(\lambda_3-\lambda_2)v}{(\lambda_1-\lambda_2)v}
\end{equation}
When $y=\cos(x)^2$, then $x=\pm\arccos(\pm\sqrt{y})$, but because $\arccos(-x)=\pi-\arccos(x)$
and the eigenvectors are in pairs $\vec{v}_i$ and $-\vec{v}_i$ and perpendicular,
the second $\pm$ can be discarded, giving:
\begin{equation}\label{phitwo}
 \phi_2 = \pm\arccos(\sqrt{v})
\end{equation}
\begin{equation}\label{phithree}
 \phi_3 = \pm\arccos(\sqrt{w})
\end{equation}
In the case that $v=0$, then $w=1$.
Because of this eigenvector symmetry, any angle $\phi_i$ is mod $\pi$,
so that when for example $\phi_2=0\mod\pi/2$, the sign of $\phi_2$ is arbitrary.
Let the following 2-dimensional vectors be defined:
\begin{equation}
 \vec{f}_1 = \left(
   \begin{matrix}
     A_{12} \\
    -A_{13}
   \end{matrix}
   \right)
\end{equation}
\begin{equation}
 \vec{f}_2 = \left(
   \begin{matrix}
     A_{22}-A_{33} \\
    -2A_{23}
   \end{matrix}
   \right)
\end{equation}
\begin{equation}
 \vec{g}_1 = \left(
   \begin{matrix}
     \frac{1}{2}(\lambda_1-\lambda_2)\cos(\phi_2)\sin(2\phi_3) \\
     \frac{1}{2}[(\lambda_1-\lambda_2)w + \lambda_2-\lambda_3]\sin(2\phi_2) 
   \end{matrix}
   \right)
\end{equation} 
\begin{equation}
 \vec{g}_2 = \left(
   \begin{matrix}
     (\lambda_1-\lambda_2)[1+(v-2)w]
      + (\lambda_2-\lambda_3)v  \\
     (\lambda_1-\lambda_2)\sin(\phi_2)\sin(2\phi_3)
   \end{matrix}
   \right)
\end{equation} 
Then the following two identities exist, where $R$ is given by (\ref{rot2}),
which are computed in section \ref{program} below by the computer algebra program:
\begin{equation}\label{phimatone}
 \vec{g}_1 = R(\phi_1)\cdot\vec{f}_1
\end{equation}
\begin{equation}\label{phimattwo}
 \vec{g}_2 = R(2\phi_1)\cdot\vec{f}_2
\end{equation}
This means that the anti-clockwise angle $\phi_1$ is equal to the
angle between $\vec{g}_1$ and $\vec{f}_1$ and half the angle between
$\vec{g}_2$ and $\vec{f}_2$. This also means that $|\vec{f}_1|=|\vec{g}_1|$
and $|\vec{f}_2|=|\vec{g}_2|$.
Let $\text{angle}(\vec{r})$ be the anti-clockwise angle of the 2-dimensional
vector $\vec{r}$ with respect to the positive $x$-axis.
Let the angles $\psi_1$ and $\psi_2$ be the following angles:
\begin{equation}
 \psi_1 = \text{angle}(\vec{f}_1)
\end{equation}
\begin{equation}
 \psi_2 = \text{angle}(\vec{f}_2)
\end{equation}
Then a rotation that rotates $\vec{f}_1$ resp. $\vec{f}_2$ onto the positive $x$-axis
makes $\phi_1$ resp. $2\phi_1$ equal to the anti-clockwise angle of $\vec{g}_1$ resp. $\vec{g}_2$.
\begin{equation}
 \phi_{1(1)} = \text{angle}(R(-\psi_1)\cdot\vec{g}_1)
\end{equation}
\begin{equation}
 \phi_{1(2)} = \frac{1}{2}\text{angle}(R(-\psi_2)\cdot\vec{g}_2)
\end{equation}
For solving the signs of (\ref{phitwo}) and (\ref{phithree}),
the sign combination with the smallest absolute difference between
$\phi_{1(1)}$ and $\phi_{1(2)}$ is selected.
Because $|\vec{f}_1|$ or $|\vec{f}_2|$ can be zero, one is selected:
\begin{equation}
 \phi_1 =
 \begin{cases}
   \phi_{1(1)} & \text{if $|\vec{f}_1|\geq|\vec{f}_2|$} \\
   \phi_{1(2)} & \text{if $|\vec{f}_1|<|\vec{f}_2|$} \\
 \end{cases}
\end{equation}
Because of equations (\ref{phimatone}) and (\ref{phimattwo}) there is at least one sign combination
of $\phi_2$ and $\phi_3$ for which $\phi_{1(1)}=\phi_{1(2)}$.
When $|\vec{f}_1|=|\vec{f}_2|=0$, the matrix is already diagonal
with two equal diagonal elements, which means that at least two of the $\lambda_i$ are equal,
which is a special case treated below.
\begin{figure}[t]
\setlength{\unitlength}{0.5cm}
\begin{picture}(25,12)(0.0,0.0)
\put(6,1){\line(0,1){8}}
\put(2,5){\line(1,0){8}}
\put(9,4.5){$\vec{g}_{11}$}
\put(6.2,8.5){$\vec{g}_{12}$}
\put(9,8){$++$}
\put(9,2){$-+$}
\put(2,8){$+-$}
\put(2,2){$--$}
\thicklines
\put(6,5){\vector(2,1){3}}
\put(6,5){\vector(2,-1){3}}
\put(6,5){\vector(-2,1){3}}
\put(6,5){\vector(-2,-1){3}}
\thinlines
\put(18,1){\line(0,1){8}}
\put(14,5){\line(1,0){8}}
\put(21,4.5){$\vec{g}_{21}$}
\put(18.2,8.5){$\vec{g}_{22}$}
\put(21,8){$++$}
\put(21,7){$--$}
\put(21,3){$+-$}
\put(21,2){$-+$}
\put(14,8){$++$}
\put(14,7){$--$}
\put(14,3){$+-$}
\put(14,2){$-+$}
\thicklines
\put(18,5){\vector(1,2){1.5}}
\put(18,5){\vector(1,-2){1.5}}
\thinlines
\end{picture}
\caption{The effect of the signs of $\phi_2$ and $\phi_3$ on $\vec{g}_1$ and $\vec{g}_2$.}
\label{figsigns}
\end{figure}
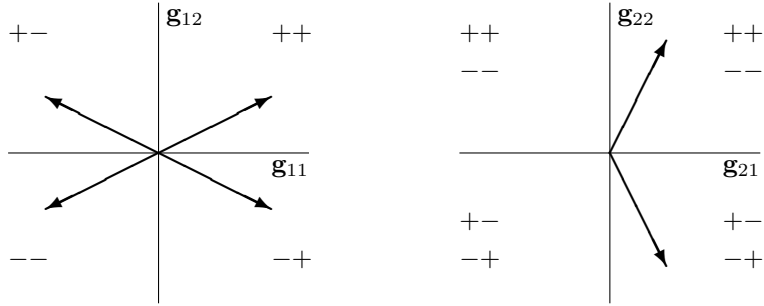

It is conjectured that (1) when $|\vec{f}_1|\neq 0$ and $|\vec{f}_2|\neq 0$, then there is
only one sign combination for which $\phi_{1(1)}=\phi_{1(2)}$,
and (2) when $|\vec{f}_1|\neq 0$ and $|\vec{f}_2|=0$, then $\phi_{1(1)}\mod\pi$ is correct for all sign combinations,
and when $|\vec{f}_1|=0$ and $|\vec{f}_2|\neq 0$, then $\phi_{1(2)}\mod\pi$ is correct for all sign combinations.
For conjecture (1), when investigating the effect of the signs of $\phi_2$ and $\phi_3$
on $\vec{g_1}$ and $\vec{g}_2$ (see figure \ref{figsigns}),
then it is clear that there may be non-unique occurences of $\phi_{1(1)}=\phi_{1(2)}$ when
(for simplicity taking $\psi_1=\psi_2=0$):
\begin{equation}
 \arctan(\frac{\vec{g}_{12}}{\vec{g}_{11}}) = \frac{1}{2}\arctan(\frac{\vec{g}_{22}}{\vec{g}_{21}})
\end{equation}
The following identity follows from $\tan(2x)=2\tan(x)/(1-\tan(x)^2)$ \cite{OLBC10}:
\begin{equation}
 \arctan(x) = \frac{1}{2}\arctan(\frac{2x}{1-x^2})
\end{equation}
which yields that in this case:
\begin{equation}
 2\vec{g}_{11}\vec{g}_{12}\vec{g}_{21} = \vec{g}_{22} ( \vec{g}_{11}^2 - \vec{g}_{12}^2 ) 
\end{equation}
This identity is computed by the computer algebra program in section \ref{program} below to be:
\begin{equation}
 (\lambda_1-\lambda_2)(\lambda_1-\lambda_3)(\lambda_2-\lambda_3) \cos(\phi_2)^2\sin(\phi_2)\sin(2\phi_3) = 0
\end{equation}
The case of at least two equal $\lambda_i$ is a special case treated below. For all $\lambda_i$ unequal
this equation means that $\phi_2=0\mod\pi/2$ or $\phi_3=0\mod\pi/2$,
which means as mentioned above that the sign of $\phi_2$ or $\phi_3$ is arbitrary.
Substituting these values in $\vec{g}_1$ and $\vec{g}_2$ shows that in these cases there are no non-unique occurences
of $\phi_{1(1)}=\phi_{1(2)}$.
For an example of conjecture (2), when $|\vec{f}_1|=0$, $A_{12}=A_{13}=0$, which means that there is a 2-dimensional
rotation in the $\vec{e}_2$-$\vec{e}_3$ plane. Let therefore $\phi_2=\phi_3=0$.
Then from the 2-dimensional matrix $D\cdot\Lambda\cdot D^T$ (\ref{mat2})
it follows that $A_{22}-A_{33}=(\lambda_2-\lambda_3)\cos(2\phi)$
and $-2A_{23}=-(\lambda_2-\lambda_3)\sin(2\phi)$.
Assuming that $\lambda_2>\lambda_3$, this means that $\psi_2=\text{angle}(\vec{f}_2)=-2\phi$.
Because in this case $\vec{g}_{22}=0$, and as $v=w=1$, $\vec{g}_{21}=\lambda_2-\lambda_3$,
which means that $\vec{g}_2$ is on the positive x-axis.
Because in this case $R(-\psi_2)=R(2\phi)$ it is then clear that $\phi_{1(2)}=\phi\mod\pi$.\\
The fourfold rotation symmetry of the eigenvectors is generated by adding $\pi$ to the angles
$\phi_2$ and $\phi_3$. Changing the order of $\lambda_1$, $\lambda_2$ and $\lambda_3$
is compensated by change of the resulting angles $\phi_1$, $\phi_2$ and $\phi_3$,
thus always generating the same set of eigenvectors.\\
A special case occurs when $p=0$, in which case $q=0$ and $\lambda_1=\lambda_2=\lambda_3$.
Then the angles can be taken $\phi_1=\phi_2=\phi_3=0$.\\
Another special case occurs when $q^2=4p^3\neq 0$, in which case two of the three $\lambda_i$ are equal.
In this case, choosing the order of the $\lambda_i$ such that:
\begin{equation}
 \lambda_1 = \lambda_2 = \lambda
\end{equation}
yields the following angles:
\begin{equation}
 s = \cos(\phi_2)^2 = \frac{A_{11}-\lambda_3}{\lambda-\lambda_3}
\end{equation}
\begin{equation}
 \phi_2 = \pm\arccos(\sqrt{s})
\end{equation}
\begin{equation}
 \phi_3 = 0
\end{equation}
This identity is computed in the computer algebra program in section \ref{program} below.
In this special case the vectors $\vec{g}_1$ and $\vec{g}_2$ simplify to:
\begin{equation}
 \vec{g}_1 = \left(
   \begin{matrix}
     0 \\
     \frac{1}{2}(\lambda-\lambda_3)\sin(2\phi_2) 
   \end{matrix}
   \right)
\end{equation} 
\begin{equation}
 \vec{g}_2 = \left(
   \begin{matrix}
     (\lambda-\lambda_3)s  \\
     0
   \end{matrix}
   \right)
\end{equation} 
The other equations are identical to the general case, except that only two sign combinations
of $\phi_2$ are possible.\\

\section{Euler Angles}\label{secteuler}

As the three angles $\phi_i$ are now determined, the Euler angles of the eigenvectors
with respect to the basis vectors can now be computed.
A 3-dimensional rotation can always be written as a sequence of three rotations
around the basis vectors as fixed reference axes as in (\ref{rotmat3}):
\begin{equation}
 D = R_1(\phi_1) \cdot R_2(\phi_2) \cdot R_3(\phi_3) 
\end{equation}
When the reference axes are also rotated, a so called Euler rotation sequence results.
When $S$ gives the rotation of the reference axes, then the new rotation $R'$
around the new reference axis is:
\begin{equation}
 R' = S\cdot R\cdot S^{-1}
\end{equation}
Thus the second rotation is a rotation around the new rotated second reference axis:
\begin{equation}
 R_2'(\phi_2) = R_1(\phi_1)\cdot R_2(\phi_2)\cdot R_1(\phi_1)^{-1}
\end{equation}
and the third rotation is a rotation around the newest rotated third reference axis:
\begin{equation}
 R_3''(\phi_3) = (R_2'(\phi_2)\cdot R_1(\phi_1))\cdot R_3(\phi_3)\cdot (R_2'(\phi_2)\cdot R_1(\phi_1))^{-1}
\end{equation}
which results in the total Euler rotation sequence:
\begin{equation}\label{eulerrot}
 R_3''(\phi_3)\cdot R_2'(\phi_2)\cdot R_1(\phi_1) = R_1(\phi_1) \cdot R_2(\phi_2) \cdot R_3(\phi_3)
\end{equation}
This equation means that a rotation sequence around fixed reference axes is equal
to an Euler rotation sequence around rotating reference axes with identical angles
but in reverse order.
The right side of this equation is the total rotation $D$ in equation \ref{rotmat3}, and
as the $R$-rotations were anti-clockwise, the anti-clockwise Euler rotation
angles thus are $\{\phi_1,\phi_2,\phi_3\}$.

\section{Computer Algebra Program}\label{program}

The Mathematica$^{\textregistered}$ \cite{W03} program used to compute the expressions
is given below.

\begin{alltt}
L:=\{\{lambda1,0,0\},\{0,lambda2,0\},\{0,0,lambda3\}\}
R1:=\{\{1,0,0\},\{0,Cos[phi1],-Sin[phi1]\},\{0,Sin[phi1],Cos[phi1]\}\}
R2:=\{\{Cos[phi2],0,Sin[phi2]\},\{0,1,0\},\{-Sin[phi2],0,Cos[phi2]\}\}
R3:=\{\{Cos[phi3],-Sin[phi3],0\},\{Sin[phi3],Cos[phi3],0\},\{0,0,1\}\}
Rtot:=R1.R2.R3
A:=Rtot.L.Transpose[Rtot]
FullSimplify[A[[1,2]]^2+A[[1,3]]^2
 +(A[[1,1]]-lambda3)(A[[1,1]]+lambda3-lambda1-lambda2)]
(lambda1-lambda3)(-lambda2+lambda3)Cos[phi2]^2
FullSimplify[A[[1,1]]-lambda3+(lambda3-lambda2)Cos[phi2]^2]
(lambda1-lambda2)Cos[phi2]^2Cos[phi3]^2
f1:=\{A[[1,2]],-A[[1,3]]\}
f2:=\{A[[2,2]]-A[[3,3]],-2A[[2,3]]\}
g1:=\{1/2(lambda1-lambda2)Cos[phi2]Sin[2phi3],
 1/2((lambda1-lambda2)Cos[phi3]^2+lambda2-lambda3)Sin[2phi2]\}
g2:=\{(lambda1-lambda2)(1+(Cos[phi2]^2-2)Cos[phi3]^2)
 +(lambda2-lambda3)Cos[phi2]^2,
 (lambda1-lambda2)Sin[phi2]Sin[2phi3]\}
R[phi]:=\{\{Cos[phi],-Sin[phi]\},\{Sin[phi],Cos[phi]\}\}
FullSimplify[R[phi1].f1-g1]
\{0,0\}
FullSimplify[R[2phi1].f2-g2]
\{0,0\}
FullSimplify[2g1[[1]]g1[[2]]g2[[1]]-g2[[2]](g1[[1]]^2-g1[[2]]^2)]
(lamda1-lambda2)(lambda1-lambda3)(lambda2-lambda3)
 Cos[phi2]^2Sin[phi2]Sin[2phi3]
FullSimplify[ReplaceAll[A[[1,1]]-lambda3,
 \{lambda1->lambda,lambda2->lambda\}]]
(lambda-lambda3)Cos[phi2]^2
\end{alltt}

\section{C++ program}

The C++ program for testing the formulas
is in the header file diagonalize.hpp and the source file diagonalize.cpp
that have been uploaded to arXiv.

\pdfbookmark[0]{References}{}

\end{document}